\documentclass[leqno, a4paper, 12pt]{article}
\usepackage{amsmath,amsfonts,here,epsf}
\usepackage[latin1]{inputenc}
\usepackage[T1]{fontenc}
\usepackage{amsmath,amssymb}
\usepackage{ae,aecompl}

\begin{document}
\newtheorem{theo}{Theorem}[section]
\newtheorem{lemme}[theo]{Lemma}
\newtheorem{cor}[theo]{Corollary}
\newtheorem{defi}[theo]{Definition}
\newtheorem{prop}[theo]{Proposition}
\newcommand{\beq}{\begin{eqnarray}}
\newcommand{\enq}{\end{eqnarray}}
\newcommand{\be}{\begin{eqnarray*}}
\newcommand{\en}{\end{eqnarray*}}
\newcommand{\Td}{\mathbb T^d}
\newcommand{\T}{\mathbb T}
\newcommand{\Rd}{\mathbb R^d}
\newcommand{\R}{\mathbb R}
\newcommand{\Zd}{\mathbb Z^d}
\newcommand{\Z}{\mathbb Z}
\newcommand{\N}{\mathbb N}
\newcommand{\Linf}{L^{\infty}}
\newcommand{\dt}{\partial_t}
\newcommand{\Dt}{\frac{d}{dt}}
\newcommand{\Dtt}{\frac{d^2}{dt^2}}
\newcommand{\demi}{\frac{1}{2}}
\newcommand{\vf}{\varphi}
\newcommand{\epu}{_{\epsilon}}
\newcommand{\ep}{^{\epsilon}}
\newcommand{\ds}{\displaystyle}
\title{Quasi-neutral limit of the Euler-Poisson and Euler-Monge-Amp\`ere systems}
\author{G. Loeper\thanks{Research supported by a doctoral grant at Laboratoire J.A. Dieudonn\'e, Universit\'e de Nice-Sophia-Antipolis}}
\date{}
\maketitle
\thispagestyle{empty}
\bibliographystyle{plain}

\begin{abstract}
This paper studies the pressureless Euler-Poisson system and its fully non-linear counterpart,  the 
Euler-Monge-Amp\`ere system, where the fully non-linear Monge-Amp\`ere equation substitutes for the linear Poisson equation.
While the first is a model of plasma physics, the second is derived as a geometric approximation to the Euler 
incompressible equations. Using 
energy estimates, convergence of both systems to the Euler incompressible equations is proved.   
\end{abstract}
\section{Introduction}

In this paper we consider a model of a collisionless plasma where the ions are supposed to be at rest and create a neutralizing background field.  The motion of the electrons can then be described by using either the kinetic formalism or the hydrodynamic equations of conservation of mass and momentum as we do here. The self-induced electric field is the gradient of a potential that depends on the electron's density $\rho$ either through the linear Poisson equation: $\ds\Delta \phi =\frac{1}{\epsilon}(\rho-1)$, or through the fully non-linear Monge-Amp\`ere equation: $\ds\det (I+\epsilon\partial_{ij}\phi)=\rho$.  This gives the Euler-Poisson ($(EP\epu)$) system  and Euler-Monge-Amp\`ere ($(EMA\epu)$)  system. 
The non-dimensional rescaled version of both systems is the following:
\be
&&\dt\rho +\nabla\cdot (\rho v)=0,\\
&&\dt v + v\cdot \nabla v=\frac{1}{\epsilon}\nabla \phi, \\
&&\epsilon\Delta \phi =\rho-1 \,\textrm{  in the Poisson case },\\
&&\det (I+\epsilon\partial_{ij}\phi)=\rho \,\textrm {  in the Monge-Amp\`ere case. }
\en
Note that the systems are pressureless, and the only force is due to electrostatic interaction.
The energy of those systems is given by
\be
{\cal E}\ep = \demi \int \rho|v|^2 + |\nabla\phi|^2 \ dx
\en for Euler-Poisson and
 
\be
{\cal E}\ep = \demi \int \rho|v|^2 + \rho |\nabla\phi|^2 \ dx
\en for Euler-Monge-Amp\`ere.

The asymptotic we look at consists in considering large scales compared to the Debye length ($\epsilon $). At those scales the plasma appears to be electrically neutral. In this limit the plasma is expected to behave like an incompressible fluid, therefore governed by the incompressible Euler equation $(E)$.
We intend to rigorously justify those limits in the present work.

\paragraph*{Physical interpretation of the quasi-neutral limit for Euler-Poisson}
The complete model of collisionless plasma describes the behavior of two species: the ions and the electrons. However the ratio of the electron's mass and the ion's mass is of several orders of magnitude, therefore we make the assumption that the ions are at rest, and distributed over a regular grid. This assumption will imply the neutralizing background electric field (the '-1' term in the Poisson equation $\ds \epsilon\Delta\phi=\rho-1$). 
The parameter $\epsilon$ comes from the vacuum permittivity, obtained after many rescalings of the equation. The typical value of $\epsilon^2$ is between $10^{-10}$ and $10^{-5}$. We consider therefore $\epsilon$ as a small parameter, and investigate the limit $\epsilon  \to 0$ of the Euler-Poisson system.
Note that if $(v,\rho,\phi)$ is a solution to the Euler-Poisson system with $\epsilon = 1$, $(v\ep, \rho\ep, \phi\ep):=(v,\rho,\phi)(\epsilon^{-1}t, \epsilon^{-1}x)$ is a solution of the Euler-Poisson system with parameter $\epsilon$. Therefore, the limit $\epsilon \to 0$ can be interpreted as a study of the long time - large scale behavior of the system.

\paragraph*{Geometric interpretation of the quasi-neutral limit for Euler-Monge-Amp\`ere}
Whereas the Euler-Poisson system relies on a well known physical model,
the Euler-Monge-Amp\`ere system,  less famous, is a fully non-linear (but asymptotically close in the quasi-neutral regime) version of the
Euler-Poisson system;
it can be seen as a non-linear model of electrostatic interaction with the advantage of allowing finite electric field for point charges (see also \cite{BrBI} where the Born-Infeld system of electromagnetism is studied, a system that exhibits similar non-linear features). 
Apart from this interpretation, the main motivation for the study of the Euler-Monge-Amp\`ere system is the following: it appears as a 'canonical' relaxation of the geodesics on the group of
measure preserving diffeomorphisms (therefore of the Euler incompressible equation, see \cite{AK}).
This interpretation will be 
developed more accurately in the section \ref{3section-euler}  devoted to the Euler-Monge-Amp\`ere equation. 
This model was first introduced, in a discrete version, by Y. Brenier in \cite{Br3}. 
Later, a kinetic version, the Vlasov-Monge-Amp\`ere system, has been introduced and studied by  Y. Brenier and the author  
in \cite{BL}.
The present work may be seen as a further step in this study.

To see why the $(EP\epu)$ and the $(EMA\epu)$ systems should be asymptotically close in the quasi-neutral regime,
notice that if $\rho$ is close to 1 then  $\epsilon\partial_{ij}\phi$ should be small, hence
$\det (I+\epsilon\partial_{ij}\phi)=1+\epsilon \Delta\phi+ O(\epsilon^2)$ and one recovers the Poisson equation. For this reason the proof of the convergence of both systems  will be very close 
and this is why we present them altogether.

\paragraph*{Related results concerning singular perturbations}
This work is concerned with the motion of slightly compressible fluids seen as singular perturbations of the Euler incompressible equation; this field has been widely investigated using different techniques:

- Traditional analysis and geometry on the group of volume preserving diffeomorphisms (see section \ref{3section-euler})  as done in \cite{Eb} where the convergence holds in $H^s$ norm, $s$ large, for well prepared initial data, restricted to the case of barotropic fluids (i.e. when the pressure is a local function of the density $\rho$,  a case different from the one studied here.) 

- Energy estimates as done in \cite{KM} again for the case of barotropic fluids, where convergence holds in all $H^s$ norms for well prepared initial data. The result has also been extended to the non-isentropic case by \cite{MeSc}.

- Pseudo-differential energy estimates as done in \cite{Gr} 
which can be seen as a pseudo-differential generalization of \cite{KM} and where the same
convergence holds for a broader class of singular perturbations, including non-local dependence between density and pressure. 

- Modulated energy techniques for convergence of  the Vlasov-Poisson system to the so-called dissipative solutions of the Euler equation, as done 
in  \cite{Br2}, and for the Vlasov-Monge-Amp\`ere system (a kinetic version of 
the Euler-Monge-Amp\`ere system) as done in \cite{BL}.
The  convergence result obtained there holds in  weighted ``$L^2$'' norms, this method has the advantage to be valid for weak solutions, and does not 
require any smoothness of the solution. We give here more details on this result:
\paragraph*{Modulated energy technique for Vlasov-Poisson and Vlasov-Monge-Amp\`ere}
Those systems are the kinetic extensions of the $(EP\epu)$ and $(EMA\epu)$ systems. They read as follows:
\beq
&&\frac{\partial f}{\partial t}+\nabla_x\left(\xi f\right)
+\frac{1}{\epsilon}\nabla_{\xi}\left( \nabla\phi)f\right)=0\label{3vma},\\
&&\Delta \phi =\frac{1}{\epsilon}(\rho-1) \,\textrm{  in the Poisson case }\label{3pv},\\
&&\det (I+\epsilon\partial_{ij}\phi)=\rho \,\textrm {  in the Monge-Amp\`ere case }\label{3mav},\\
&&f(0,\cdot,\cdot)=f^0\label{3vmainit}.
\enq
In \cite{Br2} and \cite{BL}, the following results have been obtained:
\begin{theo}
\label{3neutre}
Let $f$ be a weak solution of (\ref{3vma}, \ref{3mav}, \ref{3vmainit}) (resp. of  (\ref{3vma}, \ref{3pv}, \ref{3vmainit})) with finite energy, let 
$(t,x)\rightarrow\bar v(t,x)$ be a smooth solution of the incompressible Euler equation (\ref{3euler}) for $t\in[0,T],$ and $p(t,x)$ the corresponding pressure, let  
\begin{eqnarray*}
H_\epsilon(t)=\frac{1}{2}\int f(t,x,\xi)|\xi-\bar v(t,x)|^2 dx d\xi + E\epu(t),
\end{eqnarray*}
where $E\epu(t)=\epsilon^{-2} \int \rho |\nabla\phi|^2/2$ (resp. $E\epu(t)=\epsilon^{-2} \int |\nabla\phi|^2/2$ in the Poisson case), 
then 
\begin{eqnarray*}
H_{\epsilon}(t)\leq  C\exp(Ct)(H_{\epsilon}(0) + C\epsilon^2),\;\forall t \in [0,T].
\end{eqnarray*}
The constant $C$  depends only of the $W^{1,\infty}_x$ norm of 
$\{\bar v(s,.), p(s,.), \partial_t p(s,.), \nabla p(s,.)\,s\in [0,T]\}.$
\end{theo}

\paragraph*{Results}
Here we shall obtain by energy estimates a convergence to the Euler
incompressible system  in $\Linf_t H^s_x$ norm for any $s$ large enough (the minimal smoothness will be made precise). 
The convergence of both systems holds on the range of time 
on which the solution of Euler is smooth enough (roughly speaking, we will need at least $D^2v$ to be bounded in $\Linf$). Our work is based on the modulated energy techniques, restricted to the case of  monokinetic velocity profiles. Indeed, the quasi-neutral limit is much more difficult without this assumption, and is even known to be false in some cases (the two-streams instability, see \cite{CoGrGu}). We will mostly restrict ourselves to the case of  well-prepared initial data, but we will investigate briefly the case of non-prepared initial data: in that case the divergence part of the initial velocity is not assumed to be small, and we only assume that the initial fluctuations of the electronic density $\rho$ are of order $\epsilon$, so that the energy remains bounded.  
The electric field is expected to oscillate at frequency $\epsilon^{-1}$ and with amplitude $O(\epsilon^{-1})$. This case will be treated in section 2.3, and we will obtain that the divergence-free part of the velocity converges strongly to a solution of the Euler incompressible equation while its potential part stays bounded, but oscillates strongly with respect to time, and therefore converges weakly to 0.

We obtain also that both systems are closer to each other than they are close to the Euler
incompressible system; $(EP\epu)$ is thus a corrector in the convergence of $(EMA\epu)$ to $(E)$.

Although the operators that define the acceleration from the density 
are differential operators (and even fully non-linear in the second case), our proof
does not use the pseudo-differential formalism.
Actually, we were not able to use  the general theorem obtained by Grenier \cite{Gr} for singular perturbations: it might be because of the absence of pressure that changes the symmetrizers of the system.
After a convenient change of variable however, the system appears under a form which is strongly reminiscent (at least for the highest order terms) of the rapidly rotating fluids. This limit has been treated in \cite{Gr}.

Finally we also mention the work of Cordier  \& Grenier \cite{CG}, and Wang \cite{Ws}, where the quasi-neutral limit 'with pressure' is treated. The techniques used there do not apply here, and it is worth noting that the scaling obtained are not the same. The reader can also refer to \cite{Guo} and \cite{Pe} where different regimes of the Euler-Poisson system are studied. 

We split the rest of the paper in two sections: the first one devoted to the study of the Euler-Poisson system, and the second devoted to the study of the Euler-Monge-Amp\`ere system.

\subsection{Notation}
Hereafter $x\in \Td=\Rd/\Zd$ and $t\in \R^+$; $v(t,x)\in \Rd$ stands for the velocity and $\rho(t,x)\in \R^+$ is the macroscopic density of electrons; $\phi(t,x)\in \R$ is the electrostatic potential; $d=2$ or $3$.

It is always assumed that $\rho(t,\cdot)$ has total mass equal to 1.

The divergence of a vector field $v$ will be denoted by ${\rm div } v$ or $\nabla\cdot v$; its rotational (or curl) will be denoted by $\nabla\times v$ or ${\rm curl } v$.

The components of a vector will be denoted with superindices, i.e. $v\in \Rd=(v^i)_{i=1..d}$.

In all the paper, $[\cdot]$ will denote the integer part.

We denote respectively by $(E)$, $(EP\epu)$, $(EMA\epu)$ the systems Euler incompressible, Euler-Poisson, Euler-Monge-Amp\`ere.

\section{The Euler-Poisson system}
We consider the following Euler-Poisson system denoted by $(EP\epu)$:
\beq
&&\dt \rho +\nabla\cdot (\rho v)=0\label{3continuite},\\
&&\dt v + v\cdot \nabla v = \frac{\nabla \phi}{\epsilon } \label{3euler-poisson},   \\
&&\epsilon\Delta \phi =\rho -1\label{3poisson},
\enq
and consider the limit $\epsilon $ going to 0. 
We recall also the incompressible Euler equation $(E)$:
\beq
&&\dt v + v\cdot \nabla v = \nabla p,\nonumber   \\
&&\nabla\cdot v=0. \label{3euler}
\enq
We recall (see  \cite{Che} for a reference on the topic) that in the periodic case, the Cauchy problem for (\ref{3euler}) is well posed in $H^s(\Td)$ if $s>d/2+1$. More precisely, if $d=2$, for any divergence-free initial datum in $H^s$, there exists a unique global solution in $\Linf_{loc}(\R, H^s(\Td))$; if $d=3$ one can only prove the existence of a smooth solution in finite time, belonging to $\Linf_{loc}([0,T[,H^s(\Td))$ for some $T>0$. 
We will then prove the following:
\subsection{Result}
\begin{theo}\label{3eulerpoisson}
Let $s\in \N$ with $s\geq [d/2]+2$. 
Let $\bar v_0$ be a divergence-free vector field on $\Td$.
Let $(\bar v,p)$ be a smooth solution of the Euler incompressible system (\ref{3euler}) on $[0,T]\times \Td$,
with initial data $\bar v_0$, satisfying $\bar v\in \Linf([0,T], H^{s+1}(\Td))$. 
Let $ (\rho_0\ep, v_0\ep) $ be a sequence of initial data such that 
$\int_{\Td}\rho\ep(x)dx =1$, and such that
$\ds \left(\epsilon^{-1}(v_0\ep-\bar v_0),\epsilon^{-2}(\rho_0\ep-1)\right)$ 
is bounded in $H^s\times H^{s-1}(\Td)$. Then there exists a sequence $ (v\ep, \rho\ep)$ of 
solutions to $(EP\epu)$  with initial data $(v_0\ep, \rho_0\ep)$, belonging to  $\Linf([0,T_{\epsilon}],H^s\times H^{s-1}(\Td))$, with $\liminf_{\epsilon\rightarrow 0}T_{\epsilon}\geq T$. Moreover for any $T'<T$ and $\epsilon$ 
small enough, 
$\ds\left(\epsilon^{-1}(v\ep-\bar v),\epsilon^{-2}(\rho\ep-1)\right)$
is bounded in $\Linf([0,T'], H^{s}\times H^{s-1}(\Td))$. 
Finally when $T=+\infty$, $T\epu$ goes to infinity.
\end{theo}

{\it Remark.} 
The models that we consider are valid in domains without boundaries, and although stated in the space periodic case, we believe that our results hold true, with some technical adaptation, in  the case of the whole space.
\subsection{Proof of Theorem \ref{3eulerpoisson}}

\subsubsection{Heuristics.}
Let us introduce $(\bar v, p)$ the solution of the Euler incompressible system (\ref{3euler}) and corresponding pressure.
Note that by taking the divergence of (\ref{3euler}) the pressure is given by the following:
\be
&&\Delta p =\sum_{i,j=1}^d\partial_i \bar v ^j \partial_j\bar v ^i.
\en
We will all along the paper use the following  notation: for two vector fields $u,v$,
\be
&& \nabla  u : \nabla  v=\sum_{i,j=1}^d\partial_i  u ^j \partial_j  v ^i. 
\en
If $v$ is solution to $(EP\epu)$, we introduce also 
\be
&&v=\bar v + \epsilon v_1,\\
&&\rho=1+\epsilon^2 \rho_1.\\
\en
We suppose for now that $d=2$,
and we  take the curl and the divergence of the momentum equation, this yields
\be
&& \dt {\rm curl} v_1 + v\cdot\nabla {\rm curl} v_1 = R_1,\\
&&\dt {\rm div} v_1 + v\cdot\nabla {\rm div} v_1 = \frac{\rho_1}{\epsilon} + R_2,\\
&& \dt \rho_1 + v\cdot\nabla\rho_1 = -\frac{{\rm div}v_1}{\epsilon} + R_3.
\en
If we assume for now that $v_1, \rho_1$ and their spatial derivatives remain bounded (we do not specify in what sense yet), $R_i,i=1,2,3$ are bounded terms.
Hence the vorticity ${\rm curl}v$ is not affected by the  electric field, and 
the vector ${\bf u} = ({\rm curl}v_1, {\rm div}v_1, \rho_1)$ evolves through
\be
\dt {\bf u} + v\cdot\nabla{\bf u}= \frac{{\bf k}}{\epsilon}\times {\bf u} + {\bf R} ,
\en with ${\bf R}$ bounded, and ${\bf k} = (1,0,0)$. Under this form,   the system looks like a rapidly rotating fluid (up to the remainder ${\bf R}$, and the fact that $v\neq {\bf u}$) and the singular term
$\frac{{\bf k}}{\epsilon}\times {\bf u}$ induces time oscillations of frequency $\epsilon^{-1}$, but does not increase the energy of the perturbation, allowing energy estimates, as long as the remainder ${\bf R}$ is under control.

\subsubsection{Reformulation of the system with new unknowns}
For a vector $u\in \Rd$, we denote $u^1,...,u^d$ its components.
We define the new unknowns $\omega_1,\beta_1,\rho_1$ as follows:
\be
&& \nabla \cdot v = \epsilon \beta_1,\\
&& \rho=1+ \epsilon^2 \rho_1=1 + \epsilon \Delta\phi,\\
&& \textrm{curl } v =\omega=\bar \omega+\epsilon \omega_1,
\en
with $(\bar v, p)$ the solution of $(E)$, and $\bar \omega= \textrm{curl } \bar v $. 

We will use the following observation: 
\begin{lemme}\label{3momentep}
Let $(\rho,v)$ be a solution to $(EP\epu)$. Then, the total momentum $\int \rho(t,x)v(t,x) \ dx$ does not depend on time.
\end{lemme}

{\bf Proof.} The momentum equation can be rewritten
\be
\dt(\rho v) + \nabla\cdot(\rho v\otimes v) = \rho \nabla\phi.
\en
We just have to show that the integral of the right hand side vanishes.
For this,  we use the identity
$$(1+\epsilon^2\Delta \phi)\nabla \phi=\epsilon^2\left[\nabla\cdot (\nabla \phi \otimes \nabla \phi)-\demi\nabla|\nabla \phi|^2\right]  + \nabla \phi,$$
that yields 0 when integrated over $\Td$.

$\hfill\Box$

We note also 
\be
v=\bar v + \epsilon v_1,
\en
but $v_1$ is not really an unknown since it can be obtained from the knowledge of $\omega_1, \beta_1$, and $\epsilon\int \rho v_1 = \int \rho_0 v_0 - \int \rho \bar v$: 
Indeed when $d=2$ we have 
\beq
&&\partial_1 \beta_1 +\partial_2\omega_1 =\Delta v_1^1\label{3alphabetav1},\\
&&\partial_2 \beta_1 - \partial_1 \omega_1 =\Delta v_1^2\label{3alphabetav2}.
\enq
In the 3 dimensional case we have equations (\ref{3alphabetav1}, \ref{3alphabetav2}) replaced by
\be
(\nabla\times \nabla\times v ) +\nabla (\nabla\cdot v)=\Delta v,
\en
and thus
\beq
(\nabla\times \omega_1 )^i+\partial_i\beta_1 = \Delta v^i_1\label{3alphabetav3}.
\enq 
Note that when $d=2$ the vorticity is scalar and when $d\neq 2$ it is a vector field of $\Td$.
We show now that $v_1$ can be estimated from $\omega_1$, $\beta_1$, and $\int \rho v$.
\begin{lemme}\label{3retrieve}
Let $(\rho,v)$ be the solution at time $t$ of $(EP\epu)$ with initial datum $(\rho_0, v_0)$, let $\bar v$ be a solution at time $t$  of $(E)$ with initial datum $\bar v_0$.
Then we have, for $s\in \R$, 
\be
&&\Big\|v-\bar v - \int [v-\bar v]\Big\|_{H^{s+1}(\Td)} \leq  C \Big( \|\nabla\cdot v\|_{H^s(\Td)} +  \|\nabla\times(v-\bar v)\|_{H^s(\Td)} \Big),\\
&&\left|\int [v- \bar v]\right| \leq \\
&&C \Big( \|\rho-1\|_{L^2(\Td)} (\|\bar v\|_{L^2(\Td)}+\|\nabla\cdot v\|_{H^{-1}(\Td)} + \|\nabla\times (v-\bar v)\|_{H^{-1}(\Td)}) + \left|\int [\rho_0v_0 - \bar v_0]\right| \Big).
\en
\end{lemme}

 {\bf Proof.} Let $u_1$ be the unique vector field with zero average such that 
\be
&& \nabla\times u_1 = \nabla\times (v-\bar v),\\
&& \nabla\cdot u_1 = \nabla\cdot v.
\en
We have directly from identities (\ref{3alphabetav1}, \ref{3alphabetav2}, \ref{3alphabetav3}) that $\|u_1\|_{H^{s+1}} \leq C \Big( \|\nabla\cdot v\|_{H^s} +  \|\nabla\times (v-\bar v)\|_{H^s} \Big)$, which is the first point of the lemma.
The difference $v-\bar v -u_1=w_1$ is  a constant vector field, and we have
\be
\int \rho v& =& \int \rho \bar v + \int \rho(v-\bar v)\\ 
&=& \int \bar v + \int (\rho - 1)\bar v + \int (\rho-1) u_1 + w_1,
\en
which yields
\be
w_1 = \int [\rho_0v_0 -\bar v_0] - \int (\rho - 1)\bar v - \int (\rho-1) u_1,
\en
and the result follows.

$\hfill \Box$

We immediately deduce the following corollary:
\begin{cor}\label{3estime}
Let $(\rho, v)$ be a solution to $(EP\epu)$, let $\bar v$ be a solution to $(E)$. Then, for any $s\geq 0$,
\be
\|v-\bar v\|_{H^s} &\leq& \left|\int [\rho_0 v_0 - \bar v_0]\right|  \\
&+& C(1+ \|\rho-1\|_{H^s}) \Big( \|\nabla\cdot v\|_{H^{s-1}} + \|\nabla\times (v-\bar v)\|_{H^{s-1}}\Big) \\
&+& C\|\rho-1\|_{H^s}\|\bar v\|_{H^s}.
\en
\end{cor}

Taking the curl of equation (\ref{3euler-poisson}) we recall the following identities:
\beq
&& \dt \omega + (v\cdot \nabla) \omega  + (\nabla\cdot  v)  \omega  =0 \textrm{ when } d=2,\label{3dtalpha2}\\
&& \dt \omega + (v\cdot \nabla) \omega  + (\nabla\cdot  v ) \omega + (\omega \cdot \nabla) v=0 
\textrm{ when } d=3.\label{3dtalpha3}
\enq
When $d=2$ the $(EP\epu)$ system  then becomes:
\beq
&&\dt (\bar \omega+\epsilon \omega_1)+ v\cdot \nabla (\bar \omega+\epsilon \omega_1)
=-(\bar \omega+\epsilon \omega_1)\epsilon \beta\label{3e1},\\
&& \dt \epsilon\beta_1 + v\cdot\nabla \epsilon\beta_1 + 2\epsilon\partial_i \bar v^j \partial_j v_1^i +\epsilon^2\partial_i v_1^j \partial_j v_1^i=\frac{\Delta \phi}{\epsilon }-\partial_i \bar v^j \partial_j \bar v^i,\\
&& \dt \epsilon^2\rho_1 + v\cdot \nabla\epsilon^2\rho_1 = -(1+\epsilon^2\rho_1)\epsilon \beta_1,
\enq
whereas when $d=3$ one would have to replace equation (\ref{3e1}) by
\beq
\dt (\bar \omega+\epsilon \omega_1)+ v\cdot \nabla (\bar \omega+\epsilon \omega_1)
=-(\bar \omega+\epsilon \omega_1)\epsilon \beta
-\omega\cdot \nabla v\label{3e2}.
\enq
Noticing that $\Delta\phi=\epsilon\rho_1$, if we set 
\beq 
\tilde \rho_1 =\rho_1 -\Delta p,
\enq
we get the following system for $d=2$:
\beq
\left\{ \begin{array}{lll} 
\dt \omega_1+ v\cdot \nabla \omega_1=-\bar \omega \beta_1 - \epsilon \omega_1 \beta_1
-v_1\cdot \nabla \bar \omega,\\
\dt \beta_1+ v\cdot \nabla \beta_1=\ds{\frac{\tilde \rho_1}{\epsilon}}-2 \nabla  v_1: \nabla \bar v - \epsilon \nabla  v_1:\nabla v_1,\\
\dt \tilde \rho_1+ v\cdot \nabla \tilde \rho_1=-\ds{\frac{\beta_1}{\epsilon}}-\epsilon(\tilde \rho_1+\Delta p) \beta_1-(\dt \Delta p + v \cdot \nabla \Delta p).
\end{array}\right .\label{3newsystem}
\enq
When $d=3$, the first equation should be replaced by
\be
&&\dt \omega_1+ v\cdot \nabla \omega_1=\\
&&-\bar \omega \beta_1 -v_1\cdot \nabla \bar \omega
-\bar\omega\cdot \nabla v_1 - \omega \cdot \nabla \bar v -
\epsilon \omega_1 \beta_1 - \epsilon \omega_1\cdot \nabla v_1. 
\en

\subsubsection{Energy estimates}
We handle the energy estimates when $d=2$ but
the same result would hold when $d=3$, just with more terms.
The system can be written in the following way:
\beq
&&\dt {\bf u}\ep + \sum_i v^i \partial_i{\bf u}\ep + R\ep {\bf u}\ep = S\ep({\bf u}\ep)\label{3beau},\\
&& {\bf u}\ep(0)={\bf u}\ep_0\label{3initiale},
\enq
where $v$ is still the velocity, 
and where 
\be
{\bf u}\ep=\left(\begin{array}{c}\omega_1\\\beta_1\\ \tilde \rho_1 \end{array}\right)
, \  R\ep=\left(\begin{array}{ccc}0&0&0\\0&0&-\frac{1}{\epsilon}\\0&\frac{1}{\epsilon}&0\end{array}\right).
\en
The 'source' term $S\ep$ is given by
\be
S\ep=\left(\begin{array}{c} -\bar \omega \beta_1 - \epsilon \omega_1 \beta_1 -v_1\cdot \nabla\bar \omega\\ -2 \nabla v_1:\nabla \bar v -\epsilon\nabla v_1:\nabla v_1 \\ -\epsilon(\tilde \rho_1+\Delta p) \beta_1-(\dt \Delta p + v \cdot \nabla \Delta p)\end{array} \right).
\en
We apply $\partial^\gamma$  to equation (\ref{3beau}), where $\gamma=(\gamma^1,...,\gamma^d)$, and
$\partial^\gamma$ stands for $\frac{\partial^{|\gamma|}}{\partial x_1^{\gamma^1}...\partial x_d^{\gamma^d}}$, with $|\gamma|=\sum_{i=1}^d \gamma^i$. We get
\be
\partial_t \partial^\gamma {\bf u}\ep + v^i \partial_i \partial^\gamma {\bf u}\ep + \Sigma\ep + R\ep \partial^\gamma {\bf u}\ep =\partial^\gamma S\ep, 
\en
where
\be
\Sigma\ep=\sum_{i=1}^d \sum_{|\mu| \geq 1, \gamma \geq \mu}\partial^\mu v^i \partial^{\gamma-\mu}\partial_i {\bf u}\ep.
\en
Then we have 
\begin{lemme}\label{3bound}
If $|\gamma| >d/2$, then for $\Sigma\ep$ and $S\ep$ defined as above we have:
\be
\|\Sigma\ep(t,\cdot)\|_{L^2(\Td)}\leq C(1+\|{\bf u}\ep(t,\cdot)\|_{H^{|\gamma|}(\Td)}+\epsilon\|{\bf u}\ep(t,\cdot)\|^2_{H^{|\gamma|}(\Td)})
\en
and 
\be
\|\partial^{\gamma}S\ep(t,\cdot)\|_{L^2(\Td)}\leq C(1+\|{\bf u}\ep(t,\cdot)\|_{H^{|\gamma|}(\Td)} +\epsilon \|{\bf u}\ep(t,\cdot)\|_{H^{|\gamma|}(\Td)}^2).
\en

\end{lemme}

{\bf Proof.}

- Point 1: the proof is a straightforward adaptation of the proof of \cite{AG} p.151.
It is based on the following estimate: \cite{AG} Proposition 2.1.2 p. 100:

\begin{prop}\label{3AG}
If $u,v\in \Linf\cap H^s$ $s\in \N$, then for any $\delta, \eta, |\delta|+|\eta|=s$, one has
\be
\|(\partial^\delta u)(\partial^\eta v)\|_{L^2}\leq C(\|u\|_{\Linf}\|v\|_{H^s} + \|u\|_{H^s}\|v\|_{\Linf}).
\en
\end{prop}
Applying this result to $\partial^\mu v^i \partial^{\gamma-\mu}\partial_i {\bf u}\ep$,  $|\mu|\geq 1$,
we obtain 
\be
\|\partial^\mu v^i \partial^{\gamma-\mu}\partial_i {\bf u}\ep\|_{L^2(\Td)}\leq 
C\left(\|\nabla v^i\|_{\Linf(\Td)}\|{\bf u}\ep\|_{H^{|\gamma|}(\Td)} + \|\nabla v^i\|_{H^{|\gamma|}(\Td)}\|{\bf u}\ep\|_{\Linf(\Td)}\right).
\en
We know that $\|\cdot\|_{\Linf(\Td)}\leq C\|\cdot\|_{H^{|\gamma|}(\Td)}$ if $|\gamma| >d/2$.
Thanks to (\ref{3alphabetav1}, \ref{3alphabetav2}) we have for any $s$ 
\beq
\|\nabla v_1\|_{H^s(\Td)}\leq C (\|\omega_1\|_{H^s(\Td)} +\|\beta_1\|_{H^s(\Td)})\label{3nabv},
\enq
thus $\|\nabla v_1\|_{H^{|\gamma|}(\Td)}\leq C \|{\bf u}\ep\|_{H^{|\gamma|}(\Td)}$.
Then using that $v=\bar v + \epsilon v_1$  we conclude. 

\bigskip

- Point 2: 
We also know thanks to Proposition \ref{3AG}  that for $s>d/2$  
\be
\|\nabla v_1:\nabla v_1\|_{H^s}\leq C\|\nabla v_1\|_{\Linf}\|\nabla v_1\|_{H^s}.
\en
It follows that for $s>d/2$, we have
\be
\|S\ep\|_{H^s}\leq C(1+\|{\bf u}\ep\|_{H^s} +\epsilon \|{\bf u}\ep\|_{H^s}^2)
\en
where $C$ depends on the smoothness of the solution of (\ref{3euler}).

$\hfill \Box$

\bigskip
\noindent
Having applied $\partial^\gamma$ to (\ref{3beau}), multiplied by $\partial^\gamma {\bf u}\ep$,
 and noticed that for any $w$, $(w,R\ep w)=0$, we obtain
\be
&&\dt \demi|\partial^\gamma {\bf u}\ep|^2 + \sum_{i=1}^d\partial_i (v^i \demi|\partial^\gamma {\bf u}\ep|^2)\\
&=&   \nabla\cdot v \ \demi|\partial^\gamma {\bf u}\ep|^2 + (\partial_\gamma S\ep-\Sigma\ep)\partial^\gamma {\bf u}\ep.
\en
Since $\|\nabla\cdot v\|_{\Linf}=\epsilon \|\beta_1\|_{\Linf}\leq C\epsilon\|{\bf u}\ep\|_{H^\gamma}$ if $\gamma>d/2$,  using Lemma \ref{3bound}, and integrating over $\Td$ we  have, for any $|\gamma| \leq s$, $s>d/2$:
\be
\Dt \|\partial_\gamma {\bf u}\ep(t,\cdot)\|_{L^2(\Td)}^2 \leq C\left(1 + \|{\bf u}\ep(t,\cdot)\|_{H^{s}(\Td)}^2 + \epsilon \|{\bf u}\ep(t,\cdot)\|_{H^{s}(\Td)}^3 \right).
\en
By summing this over all multi-indexes $|\gamma| \leq s$, we can conclude using a standard Gronwall's lemma, that if the  solution $\bar v,p$ of Euler is smooth (see below for the smoothness required) on the time interval $[0,T]$, for any $T'<T$
there exists $\epsilon_0$ such that the sequence $({\bf u}\ep)_{\epsilon<\epsilon_0}$ is bounded in $\Linf([0,T'], H^s(\R^2))$. Then we have 
\be
&&{\bf u}\ep=({\rm curl} v_1, {\rm div}v_1, \rho_1 - \Delta p),\\
&&v=\bar v + \epsilon v_1,\\
&&\rho=1 +\epsilon^2 \rho_1.
\en
We use Lemma \ref{3retrieve} to get $v$ from ${\bf u}, \bar v$, and from Corollary \ref{3estime}, 
 the bound obtained on ${\bf u}$ implies a bound on
$\epsilon^{-1}(v-\bar v), \epsilon^{-2}(\rho-1).$
(Note that from the assumption on the initial data, we have 
$\left|\int \rho_0\ep v_0\ep -\bar v_0\right| \leq C\epsilon$.)

This proves Theorem \ref{3eulerpoisson}.

$\hfill \Box$

\paragraph{Minimal regularity for the limiting field}
In order to perform our computations, we need to have at least $\beta_1, \tilde\rho_1, \omega_1$ in $\Linf([0,T]\times\Td)$ and thus to have an estimate on their norms in $\Linf([0,T], H^{s}(\Td))$,
with $s > d/2$.    
Therefore we need to control $\|v_1\|_{H^{s+1}}, \|\rho_1\|_{H^s}$. 
Then we need to apply $\partial^\gamma$ to \ref{3newsystem}, with $|\gamma|=s > d/2$, 
and we need to control $\|\dt \Delta p + v\cdot\nabla \Delta p\|_{H^s}$ (this is the 'worst' term).
This implies to control $\|\bar v\|_{H^{s+2}}$.
which requires $\bar v$ to be bounded in $\Linf([0,T], H^{s+2}(\Td))$ 
with $s>d/2$. If we take $s$ integer, we ask $s\geq [d/2]+1$, and $\bar v$ must be bounded in
$\Linf([0,T], H^{[d/2]+3}(\Td))$.

{\it Remark.} Usually, modulated energy techniques only require a bound on $\|\nabla \bar v\|_{\Linf}$. Here we need one more derivative, since  the 
'div-curl' formulation of the system (performed in order to obtain energy estimates) is obtained by differentiation.

\subsection{Non-prepared initial data}

Here we obtain energy estimates in the general case of non-prepared initial data. 
What we mean by 'general case' is the case of a generic smooth initial velocity, and smooth initial density, with finite energy, hence $\rho-1$ will be of order $\epsilon$.

We will see that the energy estimates are the same as in the case of prepared initial data, the asymptotic $\epsilon \to 0$ is then handled similarly to \cite{Gr-ns}, although the algebra in our case is quite simple. The solution will exhibit a good space regularity and a strongly oscillating behavior with respect to time. 
As explained in \cite{Gr-ns}, the motion can be decomposed along a slow and a fast manifold; the slow manifold consists of divergence-free velocities with uniform density, and the fast manifold consists of potential velocities. Due to the rapid oscillations, the potential part of $v$ will converge weakly to 0, and the divergence-free part will converge strongly to a smooth solution of the incompressible Euler equation. 

We still consider $(\rho\ep, v\ep)$ solution to $(EP\epu)$.
For any vector field $v$, we introduce its soleno\"idal part and potential part, which is  the pair $(\Pi v, \nabla q)$ such that 
$v=\Pi v + \nabla q$ with $q$ periodic and $\nabla\cdot  \Pi v=0$.

\paragraph*{A priori estimates}
We first rescale the density fluctuation
\be
&&\rho\ep = 1 + \epsilon \rho_1\ep.
\en
We use the unknown ${\bf u}\ep$, given by
\be
{\bf u}\ep = \left( \begin{array}{c}\omega\ep=\text{curl} v\ep \\ \beta\ep=\text{div} v\ep  \\  \rho_1\ep\end{array} \right).
\en
(Note that we do not subtract $\Delta p $ to $\rho_1\ep$ in this case.)

We now restrict to the case $d=2$, but one can check easily that the same results hold when $d=3$.

We write the equation followed by ${\bf u}\ep$:
\beq
&&\dt {\bf u}\ep + \sum_i v^i \partial_i{\bf u}\ep + R\ep {\bf u}\ep = S\ep({\bf u}\ep)\label{3beaunp},\\
&& {\bf u}\ep(0)={\bf u}_0\ep\label{3initialenp},
\enq
where $R\ep$ is a before, and the source term $S\ep$ is now given by

\be
S\ep=\left(\begin{array}{c}-\beta\ep \omega\ep  \\ -\nabla v\ep:\nabla v\ep \\  -\rho_1\ep\beta\ep \end{array} \right).
\en
When $s > d/2$, proceeding as in Lemma \ref{3bound}, we have the estimate
\be
\|S\ep\|_{H^s} \leq C_s  \|{\bf u}\ep\|_{H^s}^2.
\en

\bigskip

Arguing as in the previous case, we conclude that, given $s\geq  [d/2]+2$ and a sequence of initial data $(\rho_0\ep, v_0\ep) $, with $\int  \rho_0\ep =1$, and  such that
\be
&&\|v_0\ep\|_{H^s}\leq C, \\
&& \epsilon^{-1}\|\rho_0\ep-1\|_{H^{s-1}} \leq C,
\en
there exists a sequence $(\rho\ep, v\ep)$ of solutions to $(EP\epu)$ with initial data $(\rho_0\ep, v_0\ep)$, 
and $(\epsilon^{-1}(\rho\ep-1),v\ep)$ remains bounded in $\Linf([0,T],H^{s-1}\times H^s(\Td))$ for some $T>0$, independent of $\epsilon$.

\paragraph*{Convergence}
We use the change of variable used for rotating fluids (see \cite{Sch}), that removes the time oscillations:
Considering the pair $\tilde {\bf u}\ep=  (\tilde \beta\ep , \tilde \rho_1\ep)$ such that
\be
\tilde \beta\ep +i\tilde \rho_1\ep = e^{it/\epsilon}\Big(\beta\ep +i\rho_1\ep\Big),
\en
we have 
\be
\dt \tilde {\bf u}\ep+ v\ep\cdot \nabla \tilde {\bf u}\ep = T\ep,
\en 
with 
\be
&&\|T\ep\|_{\Linf([0,T],H^{s-1})} \leq C,\\
&& \|v\ep\|_{\Linf([0,T],H^{s})} \leq C,\\
&&\|\tilde{\bf  u}\ep\|_{\Linf([0,T],H^{s-1})} \leq C
\en
from the a priori bounds ( we still have $s\geq [d/2]+2$).
Hence, $\dt \tilde{\bf u}\ep$ is uniformly bounded at least in $L^2$, and we deduce classically that 
$(\beta\ep,\rho_1\ep)$ converges weakly to 0 in $[0,T]\times\Td$.

We assume now that $\Pi v_0\ep$, the soleno\"idal part of $v_0\ep$, converges weakly in $L^2$ to some limit $\bar v_0$, hence it converges strongly in $H^{s'}$ for $s' < s$.
We check that $\dt \omega\ep$ is bounded uniformly on $L^2$ under our assumptions. Hence from the a priori bound, $\omega\ep$ converges (if necessary passing to a subsequence, but see the remark below) in $C([0,T], H^{s'-1}(\Td))$ for all $s' < s$.

Let $\bar v\ep$ be the unique vector field with zero average such that $\nabla \times\bar v\ep = \omega\ep$. Then, $\Pi v\ep$, the soleno\"idal part of $v\ep$ is equal to $\bar v\ep + c\ep$ with $c\ep$ a constant vector field.
Since $\int \rho\ep v\ep = \int \rho_0\ep v_0\ep$, and since $\epsilon^{-1}(\rho\ep-1)$ is bounded in $H^{s}$, we have $\lim c\ep = \lim \int v_0\ep=\lim \int \Pi v_0 \ep$.
Hence $\Pi v\ep= \bar v\ep + c\ep$ converges in $C([0,T], H^{s'}(\Td))$ for all $s' < s$.

Decomposing $v\ep$ as $v\ep= \nabla q\ep + \Pi v\ep$, we have $q\ep=\Delta^{-1}\beta\ep$ that converges weakly to 0 in $L^2([0,T]\times \Td)$.
Hence, using the a priori bounds,  in the vorticity equation
\be
\dt\omega\ep+ v\ep\cdot\nabla\omega\ep = - \beta\ep \omega\ep,
\en 
we can pass to the weak limit in $L^2([0,T]\times \Td)$ and state that $\omega=\lim \omega\ep$ satisfies in ${\cal D}'([0,T]\times \Td)$,
\be
&&\dt\omega+ \Pi v\cdot\nabla\omega = 0,\\
&& \omega(0) = \omega_0,
\en
where $\Pi v$ is the limit of the soleno\"idal part of $v\ep$.
Moreover we have $\omega = \nabla\times \Pi v$, 
hence $\Pi v$ is a solution to the incompressible Euler equation, with initial data $\bar v_0$ the limit of the soleno\"idal part of $v_0\ep$.

{\it Remark 1.} From the regularity of $v_0$, the solution (in the distribution's sense) to $(E)$ with initial data $\Pi v_0$ is unique. Therefore the whole sequence $\Pi v\ep$ is converging.

{\it Remark 2.} Here we did not introduce the solution of the limit equation (incompressible Euler in this case), and chose to argue by compactness. This method looks simpler,  however, we obtain less informations concerning the 'rate of convergence' of the sequence $(\rho\ep, v\ep)$ with respect to $\epsilon$.

We gather those results in the following theorem:

\begin{theo}\label{3eulerpoissonnp}
Let $s\geq [d/2]+2$, let $(\rho\ep_0, v\ep_0)$ be a sequence of initial data, such that $\int \rho_0\ep = 1$. Assume that  $v\ep_0$ is bounded in $H^{s}(\Td)$ and that $\epsilon^{-1}(\rho\ep_0-1)$ is bounded in $H^{s-1}(\Td)$. Then 
\begin{enumerate}
\item There exists $T>0$ such that a sequence of solutions to $(EP\epu)$ with initial data $(\rho\ep_0, v\ep_0)$ exists on $[0,T]$, and the sequence $(\epsilon^{-1}(\rho\ep-1), v\ep)$ is bounded in $\Linf([0,T],H^{s-1}\times H^s(\Td))$.
\item The potential part of $v\ep$, defined above, converges weakly to 0 in $L^2([0,T]\times \Td)$.
\item If $\Pi v_0\ep$, the soleno\"idal part of $v_0\ep$ converges to some $\bar v_0$ weakly in $L^2$, then, for all $s'<s$, $\Pi v\ep$, the soleno\"idal part of $v\ep$, converges in $C([0,T], H^{s'}(\Td))$ to $\bar v$ the solution of $(E)$ with initial data $\bar v_0$.
\end{enumerate}
\end{theo}


\section{The Euler-Monge-Amp\`ere system}
We consider here the following Euler-Monge-Amp\`ere system denoted by $(EMA\epu)$:
\beq
&&\dt \rho +\nabla\cdot (\rho v)=0\label{3conti},\\
&&\dt v + v\cdot \nabla v = \frac{\nabla \psi -x}{\epsilon^2 }\label{3ema},   \\
&&\det D^2 \psi  =\rho\label{3ma}.
\enq
The last equation must be understood in the following weak sense: $\psi$ is the only (up to a constant) convex function with $\psi -|x|^2/2$ being $\Zd$ periodic such that
\beq
\forall f \in C^0(\Td),\, \int_{\Td}f(\nabla\psi)d\rho =  \int_{\Td}f(x)dx\label{3weakma}.
\enq
This definition will be made precise in Theorem \ref{3polar-per} and Definition \ref{3MaTdrho}.

\subsection{Geometric derivation of the Euler-Monge-Amp\`ere system}
This derivation has been introduced in \cite{Br3},  \cite{BL}, \cite{These}. We reproduce it for sake of completeness, but the reader only interested in the proof of  convergence may skip this section.
\label{3section-euler}
\subsubsection{The Euler equations of an incompressible perfect fluid}
The motion of an inviscid and incompressible fluid in a domain $\Omega\subset\Rd$ is described by the
 Euler incompressible equation $(E)$ that we recall here:
\be
&&\dt v +(v\cdot\nabla) v=\nabla p,\\
&&\nabla\cdot v =0.
\en
Following Arnold (see \cite{AK}), we have a formal interpretation of the Euler incompressible equations:
 introducing $G(\Omega)$ the group of all volume preserving diffeomorphisms of $\Omega$
 with jacobian determinant equal to 1, the Euler equations describe the geodesics of $G(\Omega)$
with length measured in the $L^2$ sense.

\subsubsection{Approximate geodesics }\label{3section-apxgeo}
A general strategy to define approximate geodesics along a manifold $M$
(in our case $M=G(\Omega)$) embedded in a Hilbert space $H$ 
(here $H=L^2(\Omega,\Rd)$)
is to introduce a penalty parameter $\epsilon>0$
and the following $unconstrained$ dynamical system in $H$
\beq
\partial_{tt}{X} +\frac{1}{2\epsilon^2}\nabla_X\left(d^2(X,M))\right)=0.
\label{3apxgeo}
\enq
In this equation,
the unknown $t\rightarrow X(t)$ is a curve in $H$, 
$d(X,M)$ is the distance (in $H$) of $X$ to the manifold
$M$, i.e. in our case when $M=G(\Omega)$,
\beq
d(X,G(\Omega))=\inf_{g\in G(\Omega)}\|X-g\|_H,\label{3projection}
\enq
finally, $\nabla_X$ denotes the
gradient operator in $H$.
This penalty approach 
has been introduced for the Euler equations 
by Brenier in \cite{Br3}. 
It is similar-but not identical-
to Ebin's slightly compressible flow theory \cite{Eb}, and is
a natural extension of the theory of 
constrained finite dimensional mechanical systems \cite{RU}. 
Actually if $G(\Omega)$ were a smooth manifold, the result would be
exactly the one of \cite{Eb}, Theorem 2.7, but this is not the case, 
here because the $L^2$ metric is too weak.
The penalized system is formally hamiltonian in variables
$(X,\partial_t X)$ with 
hamiltonian (or energy) given by:
\be
E=\frac{1}{2}\|\partial_t X\|_H^2 + \frac{1}{2\epsilon^2}d^2(X,G(\Omega)).
\en
Multiplying equation (\ref{3apxgeo}) by $\partial_t X$, we get immediately 
that the energy is formally conserved. 
Therefore it is plausible that the map $X(t)$ will remain close to $G(\Omega)$ 
if it is close at $t=0$.
A formal computation shows that, given a point $X$ for which there is
a unique closest point $\pi_X$ to $X$ in the $H$ closure
of $G(\Omega)$, 
we have:
\beq
\nabla_X\left(d(X,G)\right)=\frac{1}{d(X,G)}(X-\pi_X).\label{3grad}
\enq
Thus the equation (\ref{3apxgeo}) formally becomes: 
\be
\partial_{tt} X + \frac{1}{\epsilon^2}(X-\pi_X)=0\label{3apxgeo2}.
\en
To understand why solutions to such a system may approach
geodesics along $G(\Omega)$ as 
$\epsilon$ goes to 0, just recall that, in the simple framework of a 
surface $S$ embedded in the 3 dimensional Euclidean space,
a  geodesic $t\rightarrow s(t)$  along $S$ 
is characterized by the fact that
for every $t$, the plane defined by $\{\dot s(t), \ddot s(t)\}$ is 
orthogonal to $S$.
In our case, $\partial_{tt} X(t)$ is orthogonal to $G(\Omega)$ thanks to 
(\ref{3projection}) and $X(t)$ remains close to $G(\Omega)$.

\paragraph*{Notation}
Since we intend to work on the flat torus $\Td$ we might consider
$\Zd$ additive mappings, i.e. mappings such that 
$$\forall \vec{p}\in \Zd, \  X(.+\vec{p})=X(.)+\vec{p},$$ 
as well as periodic mappings (i.e. mappings from $\Td$ into itself).
\\
Then given $m$ an additive mapping,we denote by $\hat m$ the naturally associated
 mapping on $\Td$.
The following polar factorization Theorem
is a periodic version of \cite{Br1}, it has been discovered independently by \cite{Mc2} and \cite{Co}.
\begin{theo}\label{3polar-per}
Let $X: \Rd \rightarrow \Rd$ be $\Zd$ additive 
and assume  that $\rho_X= X_{\#}dx$ has a density  in  $L^1([0,1]^d)$, 
then there exits an a.e. unique pair $(\nabla\phi_X, \pi_X)$   satisfying 
$$X=\nabla\phi_X\circ \pi_X$$ 
with $\phi_X$ a convex function such that
$\phi_X(x)-|x|^2/2$ is $\Zd$ periodic, and $\pi_X: \Rd \rightarrow \Rd$  
additive such that $\hat \pi_X$ is measure preserving in $\Td$.  
Moreover we have 
$$\|X- \pi_X\|_{L^2([0,1]^d)}=\|\hat X- \hat \pi_X\|_{L^2(\Td)}=d(\hat X,G(\Td))$$ 
and if $\psi_X$ is the Legendre transform of $\phi_X$ then 
$$\pi_X=\nabla\psi_X\circ X.$$
\end{theo}

{\it Remark 1.} The pair $(\phi_X, \psi_X)$ is uniquely defined by the density 
$\rho_X=X_{\#}dx$.

{\it Remark 2.} Important properties of the optimal potential: The periodicity of $\phi_X(x)-|x|^2/2$
implies that $\nabla\phi_X$ and $\nabla\psi_X$ are $\mathbb Z^d$ additive, and that 
$\psi_X -|x|^2/2$ is also $\mathbb Z^d$ periodic.
This allows the following definition:
\begin{defi}\label{3MaTdrho}
Let $\rho$ be a probability  measure on $\Td$, then we denote $\phi[\rho]$ 
(resp. $\psi[\rho]$) the unique up to a constant convex function such that
 \be
&& \phi[\rho]-|\cdot|^2/2 \mbox{  is } \Zd \mbox{  periodic },\\ 
&&\forall f \in C^0(\Td),\;\int_{\Td} f(\hat{\nabla\phi}[\rho](x))dx=\int_{\Td} f(x)d\rho(x)
\en
(resp. its Legendre transform).
\end{defi}

{\it Remark.} We recover thus that $\psi[\rho],\phi[\rho]$ will be generalized 
solutions of the following Monge-Amp\`ere equations 
\be
&&\det D^2\psi=\rho,\\
&&\rho(\nabla\phi) \det D^2\phi =1.
\en

\subsubsection{Result}

We are now ready to state the main result of this section.

\begin{theo}\label{3th-ema}
Let $s\in \N$ with $s\geq [d/2]+ 2$. Let $\bar v_0$ be a divergence-free vector field on $\Td$, let $(\bar v,p)$ be a smooth solution of the Euler incompressible system (\ref{3euler}) on $[0,T]\times \Td$,
with initial data $\bar v_0$ and satisfying $\bar v\in \Linf([0,T], H^{s+1}(\Td))$.
Let $(v_0\ep, \rho_0\ep)$ be a sequence of initial data, with $\int \rho_0\ep=1$, and  such that
$\ds\left(\epsilon^{-1}(v_0\ep-\bar v_0), \epsilon^{-2}(\rho_0\ep-1)\right)$ 
is bounded in $H^s\times H^{s-1}(\Td)$. Then there exists a sequence $ (v\ep, \rho\ep)$ of 
solutions to $(EMA\epu)$  with initial data $(v_0\ep, \rho_0\ep)$ belonging to  $\Linf([0,T_{\epsilon}],H^s\times H^{s-1}(\Td))$ 
with $\liminf_{\epsilon\rightarrow 0}T_{\epsilon}\geq T$. Moreover for any $T'<T$ and $\epsilon$ 
small enough, 
$\ds\left(\epsilon^{-1}(v\ep-\bar v), \epsilon^{-2}(\rho\ep-1)\right)$
is bounded in $\Linf([0,T'], H^{s}\times H^{s-1}(\Td))$. 
Finally when $T=+\infty$, $T\epu$ goes to infinity.
\end{theo}

Before entering the proof of this result, we need some preliminary results concerning the linearization of the Monge-Amp\`ere operator.

\subsection{Linearization of the Monge-Amp\`ere operator in $H^s$ norm}
This section is devoted to the proof of the following Theorem:
\begin{theo}\label{3maHs}
Let $\rho$ be a probability measure  on $\Td$, $d\leq 3$,  let $\psi$ satisfy
\be
\det D^2\psi =\rho
\en
in the sense of  Definition \ref{3MaTdrho}. Then, there exists $\epsilon_0>0$ such that if $\|\rho-1\|_{H^2(\Td)} \leq \epsilon_0$,
 for any $s\in \N$ with $s>d/2$ there exists $C(s)$ that satisfies
\beq
&&\|D^2 \psi-I\|_{H^{s}(\Td)}\leq C(s)\|\rho-1\|_{H^{s}(\Td)}\label{3detdelta},\\
&&\|(\Delta\psi - d)-(\rho-1)\|_{H^s(\Td)}\leq C(s) \|(\rho-1)\|^2_{H^s(\Td)}\label{3detdelta2}.
\enq
\end{theo}

\subsubsection{Preliminary results}

We first state the following result obtained from \cite{Ca1}
on the regularity of solutions 
to Monge-Amp\`ere equation, adapted to the periodic case. This result will be the starting point of the proof of Theorem \ref{3maHs}.

\begin{theo}\label{3regper}
Let $\rho \in C^{\alpha}(\Td)$ for some $\alpha>0$, with $0<m\leq\rho\leq M$, be a probability measure on $\Td$,
let $\psi=\psi[\rho]$ in the sense of  Definition \ref{3MaTdrho}. Then $\psi$ is a classical solution of 
\be
\det D^2\psi=\rho
\en
and satisfies for any $\alpha '<\alpha$:
\begin{eqnarray} 
&&\|\nabla\psi(x)-x\|_{\Linf}\leq C(d)=\sqrt d /2,\\
&&\|D^2\psi\|_{C^{\alpha'}}\leq K(m,M,\|\rho\|_{C^{\alpha}}, \alpha, \alpha ').
\end{eqnarray}
\end{theo}
Then we state a classical result of elliptic regularity that we will need during the course of the proof. It can be found in \cite{GT}, Theorem 9.11.
\begin{theo}\label{3calzig}
Let $\Omega$ be an open set in $\Rd$, $u\in W^{2,p}_{loc}(\Omega)\cap L^p(\Omega), 1<p<\infty$, be a 
strong solution of the equation
\be
\sum_{i,j=1}^d a^{ij}\partial_{ij} u = f
\en
in $\Omega$ where the coefficients $a^{ij}$ satisfy
\be
&& a^{ij} \in C^0(\Omega), f \in L^p(\Omega);\\
&& \lambda |\xi|^2 \leq a^{ij}\xi_i\xi_j \leq \Lambda |\xi|^2 \ \ \forall \xi \in \Rd,
\en
for $i,j =1..d$, with $0<\lambda, \Lambda < \infty$. Then for any $\Omega'\subset\subset \Omega$, 
\be
\|u\|_{W^{2,p}(\Omega')}\leq C(\|u\|_{L^p(\Omega')}+\|f\|_{L^p(\Omega')}),
\en
where $C$ depends on $d,p,\lambda, \Lambda, \Omega', \Omega$ and the moduli of continuity of the 
coefficients $a^{ij}$ on $\Omega'$.
\end{theo}

\subsubsection{Proof of Theorem \ref{3maHs}}

\paragraph*{Sketch of the proof} 

We assume here $d=3$. We recall that $\psi$ satisfies $\det D^2\psi = \rho$.
We first have to prove that $\rho\in H^s$ implies $D^2\psi\in H^s$.
We will proceed by induction.
We recall first that for $A,B$ two $d\times d$ matrices, we have the following expansion:
\be
\det (A+tB) =\det A + t \  \text{trace}\left(A^t_{com} B\right) + o(t),
\en
where $A_{com}$ is the matrix whose elements are the minors of $A$, or co-matrix of $A$. 
Hence the elements of $A_{com}$ are polynomials of degree $d-1$ in the elements of $A$.
When $A$ is invertible, we have $A_{com} = \det A \ [A]^{-1}$.

Differentiating $s$ times the Monge-Amp\`ere equation, and denoting $M$ the co-matrix of $D^2\psi$, we will have
\be
\text{trace} \ \left(M D^2 \partial^s \psi\right) + {\cal T} = \partial^s \rho,
\en
where the first term contains the highest derivatives, and ${\cal T}$ will consist of products involving three derivatives of $\psi$. The order of each derivative will smaller or equal to $s-1$, and the sum of the three orders will be equal to $s$.
By a careful analysis, this product will be controlled  in $L^2$ by $\|D^2\psi\|_{H^{s-1}}$ and $\|\rho\|_{H^s}$, using Sobolev injections. 
From Theorem \ref{3regper}, assuming a minimal regularity for $\rho$ (i.e.  the bound (\ref{3maprop2})), $D^2\psi$, and therefore $M$ will be continuous elliptic matrices.
Hence $\partial^s \psi$ solves an elliptic problem, with continuous coefficients, and we will use the Theorem \ref{3calzig} to obtain $D^2\partial^s\psi\in L^2$.

This intermediate step will be done in Lemma \ref{3inter};
using this a priori estimate and the continuity method, we will obtain the estimate (\ref{3detdelta}). 

Then, the expansion $\det (I+D^2\varphi)= 1+ \Delta\varphi + P(\partial_{ij}\varphi)$, where $P$ is a polynomial
in $\partial_{ij}\varphi$ whose terms are of degree two or three (when $d=3$), 
will yield (\ref{3detdelta2}).

\paragraph*{Rigorous proof}
We recall that $\psi$ satisfies
\beq
&&\det D^2\psi =\rho\label{3maprop},\\
&& \|\rho-1\|_{H^2}\leq \epsilon_0\label{3maprop2},
\enq
for some $\epsilon_0$ to be chosen later.
We suppose $d=3$ and the proof can be reproduced in the case $d=2$ with minor modifications.
The parameter $\epsilon_0$ is chosen such that (\ref{3maprop2}) implies 
\be
\lambda_1\leq \rho \leq \lambda_2
\en
for some $\lambda_1>0, \lambda_2 > 0$.
Note also that thanks to (\ref{3maprop2}), $\rho$ is in $C^{\alpha}$ for $\alpha=\demi$.  Then from Theorem \ref{3regper}, $D^2\psi\in C^{\alpha'}$ with $\alpha'< \alpha$.
 Note also that since $\rho \in [\lambda_1,\lambda_2]$ and using equation (\ref{3maprop}), $D^2\psi\in C^{\alpha'}$ implies that $[D^2\psi]^{-1} \in C^{\alpha'}$, and thus $M^{ij}$ the co-matrix of $D^2\psi$
is uniformly elliptic and  in $C^{\alpha '}$.

\bigskip
\noindent
We first prove by induction that if $\gamma \in \N^d$ then 
$\rho \in H^{|\gamma|}$  implies $D^2\psi \in H^{|\gamma|}$. It can be checked during the proof that this bound will be uniform under the condition (\ref{3maprop2}) for $\epsilon_0$ small enough. 

\begin{lemme}\label{3inter}
Under assumption (\ref{3maprop2}), for any $\gamma \in \N^d$,  $\rho \in H^{|\gamma|}$ implies that $\partial^{\gamma}D^2\psi \in L^2$. If moreover $\rho \in W^{|\gamma|,6}$  then $\partial^{\gamma}D^2\psi \in L^6$.
\end{lemme}

{\bf Proof.} This lemma will be proved by induction. We first deal with the cases $|\gamma|=0,1,2$.

\bigskip

The case $\gamma=0$ is a consequence of Theorem \ref{3regper}.

\bigskip

For $|\gamma|=1$ we differentiate (\ref{3maprop}) with respect to $x_\nu$, to obtain
\beq\label{3d1}
M^{ij}\partial_{ij}(\partial_\nu \psi)=\partial_\nu \rho,
\enq
with $M^{ij}$ the co-matrix of $\partial_{ij}\psi$.
Then if $\partial_\nu \rho \in L^2$, by Theorem \ref{3calzig}, $\partial_\nu\psi \in W^{2,2}$.
If $\partial_\nu \rho \in L^6$ we also get that $\partial_\nu\psi \in W^{2,6}$.

\bigskip

For $|\gamma|=2$ differentiating once more with respect to $x_\beta$ we obtain
\beq\label{3d2}
M^{ij}\partial_{ij}(\partial_{\nu\beta} \psi)+
(\partial_\beta M^{ij})\partial_{ij}(\partial_\nu \psi),
=\partial_{\nu\beta} \rho
\enq
still with $M^{ij}$ the co-matrix of $\partial_{ij}\psi$.
Suppose that $\rho \in H^2$, then $W^{2,2}\subset W^{1, \frac{2d}{d-2}}=W^{1, 6}$ if $d=3$, 
and  $\partial_\nu\psi \in W^{2,6}$. The term $\partial_\beta M^{ij}$ is a sum of terms of the form $\partial_{ij}(\partial_\beta\psi) \partial_{kl}\psi$ and 
the second term of the left hand side of (\ref{3d2}) is thus bounded in $L^2$.  Then once again by Theorem \ref{3calzig} one gets that $D^2\partial_{\nu\beta} \psi\in  L^2$ if  $\partial_{\nu\beta} \rho\in L^2$.

Moreover if $D^2\rho\in L^6$ then $\partial_\nu\rho \in C^{\alpha}$ for some $\alpha>0$. Using (\ref{3d1}) and Schauder interior estimates (see \cite{GT}, Theorem 6.2.), we obtain
$D^2\partial_\nu\psi \in C^{\alpha'}$  .  Thus  $(\partial_\beta M^{ij})\partial_{ij}(\partial_\nu \psi)\in C^{\alpha'}$. 
From (\ref{3d2}) and Theorem \ref{3calzig} we obtain 
$\partial_{\nu\beta}D^2\psi \in L^6$.   

\bigskip

As we just saw, Lemma \ref{3inter} is true for $|\gamma|=0,1,2$.
We assume that it holds for all $\gamma$ with $|\gamma| \leq n$ for some $n\geq 2$. 
Take now $|\gamma|=n+1\geq 3$, $\rho \in H^{|\gamma|}$, and apply $\partial^{\gamma}$ to 
(\ref{3maprop}):
\beq\label{3dgammadet}
M^{ij} \partial_{ij}\partial^\gamma \psi + \sum_{ \begin{array}{c}\scriptstyle\gamma_1+ \gamma_2 + \gamma_3=\gamma\\
\scriptstyle|\gamma|-1 \geq |\gamma_1|\geq  |\gamma_2| \geq |\gamma_3|\end{array}}* \ \partial_{ij}\partial^{\gamma_1}\psi \partial_{kl}\partial^{\gamma_2}\psi \partial_{mp}\partial^{\gamma_3}\psi=\partial^{\gamma}\rho
\enq
with $*$ some constant coefficients. We call $\cal T$ the second term of the left hand side of (\ref{3dgammadet}).
Since $\rho \in H^{|\gamma|}$, $\rho \in W^{|\gamma|-1,6}$, and  we have
 $\partial^\alpha D^2\psi \in L^6(\Td)$ for any $|\alpha|\leq n$ using the induction hypothesis.
Therefore ${\cal T} \in L^2$
and since $\partial^{\gamma}\rho \in L^2$ we obtain  $M^{ij}\partial_{ij}\partial^\gamma \psi \in L^2$.
Using  Theorem \ref{3calzig} it follows that $\partial^{\gamma}D^2\psi \in L^2$.

Remember that $|\gamma|\geq 3$ thus $|\gamma_3|\leq \frac{1}{3}|\gamma| \leq \gamma -2$, and
$|\gamma_2|\leq \frac{1}{2}|\gamma| \leq \gamma -2$. 
Since $d=3$, we have $H^2 \subset C^\alpha$ for some $\alpha>0$ and thus  $\partial_{kl}\partial^{\gamma_2}\psi$,  $\partial_{mp}\partial^{\gamma_3}\psi$ are in $C^\alpha$, moreover $H^1 \subset L^6$ and  since $|\gamma_1|\leq |\gamma|-1$, $\partial_{ij}\partial^{\gamma_1}\psi$ is in $L^6$. Therefore $\cal T$ is in $L^6$.
If $\partial^\gamma\rho \in L^6$ we have $\partial^{\gamma}D^2\psi$ in $L^6$. 
Hence the lemma holds for all $|\gamma| \leq n+1$.
This achieves the proof of Lemma \ref{3inter}.

$\hfill \Box$

\bigskip

Now by induction on $|\gamma|$ we prove (\ref{3detdelta}) and (\ref{3detdelta2}).
From (\ref{3maprop2}) we have $\|\rho-1\|_{L^2}\leq \epsilon_0$ small. Take $\psi=|x|^2/2 + \varphi$ solution of $\det D^2\psi = \rho$ with $\varphi$ periodic and $\int_{\Td}\varphi =0$.
We begin to show that $\|\varphi\|_{C^{2,\alpha}}$ for some $\alpha>0$, is controlled by $\|\rho-1\|_{H^2}$.
Indeed, the periodic solution of 
\be
\det (I+D^2\varphi)=\rho
\en
can be built by the continuity method (see \cite{GT}).
Starting from $\rho_0 = 1, \varphi_0 = 0$, we use the implicit function Theorem to obtain the solution $\varphi_t$ of
\beq\label{3contin}
\det (I+D^2\varphi_t)=t\rho + (1-t).
\enq
For this we differentiate (\ref{3contin}) with respect to $t$, to obtain
\beq\label{3malin}
M^{ij}_t \partial_{ij}\partial_t\varphi_t = \rho -1,
\enq
for $t\in [0,1]$,
where $M_t$ is the co-matrix of $I+D^2\varphi_t$.
We know, from the a priori estimate of Theorem \ref{3regper}, that for $\rho\in C^{1/2}$, $D^2\varphi_t$ and therefore $M^{ij}_t$ are $C^\alpha$ elliptic matrices for all $\alpha<1/2$. 

To see why (\ref{3malin}) indeed admits a unique (up to a constant) periodic solution, we recall that $M$ is  the comatrix of a Hessian matrix, therefore it is 'divergence-free':
\be
\forall i\in [1..d], \ \sum_j M^{ij}=0.
\en
Hence equation (\ref{3malin}) can be rewritten in divergence form 
\be
\sum_{i,j}\partial_i(M^{ij}_t \partial_{j}\partial_t\varphi_t) = \rho -1,
\en
and the operator ${\cal L}= M^{ij}_t\partial_{ij}\cdot$ is a self adjoint operator on $H^1(\Td)$ and induces a bounded coercive bilinear form on $H^1_0(\Td)$, where the subscript 0 means that we impose the mean value to be 0. Then the existence/uniqueness of a solution to (\ref{3malin}) in $H^1_0$ follows by Lax-Milgram Theorem. 

Hence, $\partial_t \varphi_t$ is the unique (up to a constant) periodic solution of the above elliptic problem,  and from Schauder interior estimates we obtain 
\be
\|\partial_t \varphi\|_{C^{2,\alpha}} &\leq& C\|\rho -1\|_{C^{1/2}}\\
& \leq&  C\|\rho -1\|_{H^2},
\en uniformly in $t\in [0,1]$, and finally for $\alpha<1/2$, 
\be
\|\varphi\|_{C^{2,\alpha}} \leq C\|\rho -1\|_{H^2}.
\en

\bigskip

Then we have
\be
\det (I+D^2\varphi)=1 + \Delta\varphi + R_{ij}\partial_{ij}\varphi
\en
where $R$ is a symmetric matrix whose coefficients are polynomials in $\partial_{ij}\varphi$ of degree larger or equal to 1.
The norms $\|R_{ij}\|_{C^\alpha}$ are controlled by $\|\rho-1\|_{C^\alpha}\leq \epsilon_0$, hence, for $\epsilon_0$ small enough,  the matrix $\delta_{ij} + R_{ij}$ is uniformly bounded, elliptic, and $C^\alpha$ continuous.
Since $\varphi$ satisfies
\beq\label{3eqfi}
(\delta_{ij} + R_{ij})\partial_{ij}\varphi =\rho -1
\enq
it follows from Theorem \ref{3calzig} that $\|\partial_{ij}\varphi\|_{L^2} \leq C \|\rho-1\|_{L^2}$ and this proves 
(\ref{3detdelta}) for $\gamma=0$.

\bigskip

If $|\gamma|=1$, we have
\be
(M^{ij})\partial_{ij}\partial_\nu\varphi =\partial_\nu\rho 
\en 
with $M$ uniformly bounded, elliptic and $C^{\alpha}$ continuous. For the same reasons we have
$\|\partial_{ij}\partial_\nu\varphi\|_{L^2} \leq C \|\partial_\nu\rho\|_{L^2}$.

\bigskip

If $|\gamma|=2$, we do as in the Proof of Lemma \ref{3inter}: using (\ref{3d1}, \ref{3d2}) and keeping track of the bounds, we get 
\be
\|\partial_{ij}\partial_{\nu\beta}\varphi\|_{L^2} &\leq &C \|\partial_{\nu\beta}\rho\|_{L^2}+ C \|\rho-1\|_{H^2}^2\\
&\leq& C\|\rho-1\|_{H^2},
\en
with $C$ uniform under the assumption of Lemma \ref{3inter}.

\bigskip

If $|\gamma|\geq 3$, we go back to equation (\ref{3dgammadet}):  ${\cal T}$  is a sum of terms which contain all
a product of at least  two derivatives of $\psi$ of degree higher than 3.
Since $D^3\psi=D^3\varphi$, we have 
\be
\partial_{ij}\partial^{\gamma_1}\psi \partial_{kl}\partial^{\gamma_2}\psi \partial_{mp}\partial^{\gamma_3}\psi= \partial_{ij}\partial^{\gamma_1}\varphi \partial_{kl}\partial^{\gamma_2}\varphi \partial_{mp}\partial^{\gamma_3}\psi.
\en
We assume  by induction that $\|\partial_{ij}\varphi\|_{H^{|\gamma|-1}}\leq C\|\rho-1\|_{H^{|\gamma|-1}}$. Since $\rho-1 \in H^{|\gamma|}$ we also have that $D^2\psi \in H^{|\gamma|}$, with a uniform bound thanks to Lemma \ref{3inter}.
We remember that when $|\gamma|\geq 3$, we have $|\gamma_2| \leq |\gamma|-2, |\gamma_3| \leq |\gamma|-2$,  
therefore, using the injection of $H^2$ in $C^{1/2}$ when $d=3$, 
$\partial_{kl}\partial^{\gamma_2}\varphi$ and $\partial_{mp}\partial^{\gamma_3}\psi$ are uniformly bounded in $\Linf$. We obtain that
\be
\|\partial_{ij}\partial^{\gamma_1}\psi \partial_{kl}\partial^{\gamma_2}\psi \partial_{mp}\partial^{\gamma_3}\psi\|_{L^2} \leq  C\|\rho-1\|_{H^{|\gamma|-1}}. 
\en
Then $\partial^{\gamma}\varphi$ satisfies 
\be
(M^{ij})\partial_{ij}\partial^{\gamma}\varphi=\partial^\gamma \rho -{\cal T} 
\en 
thus, using Theorem \ref{3calzig},
\be
\|\partial^\gamma D^2\varphi\|_{L^2}&\leq& C(\|\rho-1\|_{H^{|\gamma|-1}}+ \|\partial^\gamma \rho\|_{L^2})\\
&\leq&C(\|\rho-1\|_{H^{|\gamma|}}),
\en
and we conclude that 
\be
\|D^2\psi-I\|_{H^s}\leq C(s)\|\rho-1\|_{H^s}
\en
for $s\in \N, s\geq 2$ and under condition (\ref{3maprop2});
thus (\ref{3detdelta}) is obtained.

\bigskip

Using Proposition \ref{3AG} and the fact that $\partial_{mp}\partial^{\gamma_3}\psi$ is  uniformly bounded in $\Linf$ thanks to Lemma \ref{3inter}, we can also  obtain that
\be
\|\partial_{ij}\partial^{\gamma_1}\varphi \partial_{kl}\partial^{\gamma_2}\varphi \partial_{mp}\partial^{\gamma_3}\psi\|_{L^2} &\leq& C \|D^2\varphi\|_{H^{|\gamma|}}^2  \\ 
&\leq& C\|\rho-1\|^2_{H^{|\gamma|}}
\en
for $|\gamma|\geq 3$. (When $|\gamma|=2$, the estimate holds also, but not using Proposition \ref{3AG}.) Therefore, for all $\gamma$, we have $\|{\cal T}\|_{L^2} \leq C \|\rho-1\|_{H^{|\gamma|}}^2$.

\bigskip

To conclude (\ref{3detdelta2}), we now write
$\det (I+D^2\varphi) =\rho$ under the form
\be
\Delta \varphi  = \rho -1 + P(\partial_{ij}\varphi),
\en
with $P$ consisting of products of two or three second derivatives of $\varphi$. Hence, under assumption (\ref{3maprop2}), using Proposition \ref{3AG}, we have for $s>d/2$, $\|P(\partial_{ij}\varphi)\|_{H^s} \leq C \|D^2\varphi\|_{H^s}^2$.
Using the bound (\ref{3detdelta}), we conclude that for $s\in \N, s> d/2$,
\be
\|\Delta \varphi - (\rho-1)\|_{H^s} \leq \|\rho-1\|^2_{H^s}.
\en
and 
Theorem \ref{3maHs} is proved.
  

$\hfill \Box$

\subsection{Energy estimates and proof of the convergence}

The proof of the energy estimates for Euler-Monge-Amp\`ere is much inspired from the proof of Theorem \ref{3eulerpoisson} for the following reason: by taking the divergence of equation (\ref{3ema}) one gets:
\be
\dt (\nabla\cdot v) + v\cdot \nabla (\nabla\cdot v) + \partial_i v^j \partial_j v^i 
=\frac{\Delta\psi -d}{\epsilon^2}.
\en
Suppose that $\rho$ is close to 1 at an order $\epsilon^2$ as is the case for Euler-Poisson, we guess (from Theorem \ref{3maHs}) that we have the following:
\be
&& \psi =|x|^2/2 + \epsilon^2 \varphi,\\
&& \rho =\det D^2\psi=1+\epsilon^2 \Delta\varphi +O(\epsilon^4), 
\en
and thus $\ds \Delta \psi =d + \epsilon^2 \Delta\varphi = d+\rho-1 + O(\epsilon^4)$. 
Therefore we expect that
\be
\dt (\nabla\cdot v) + v\cdot \nabla (\nabla\cdot v) + \partial_i v^j \partial_j v^i 
=\frac{\rho-1}{\epsilon^2} + O(\epsilon^2),
\en
and that the technique of Theorem \ref{3eulerpoisson}  will apply.

Before performing the div-curl decomposition for the energy estimates, we need to establish the analog of Lemma \ref{3momentep} in the present case, so that Lemma \ref{3retrieve} and its Corollary \ref{3estime} hold.
\begin{lemme}\label{3momentema}
Let $(\rho,v)$ be a solution to $(EMA\epu)$. Then, the total momentum $\int \rho(t,x)v(t,x) \ dx$ does not depend on time.
\end{lemme}

{\bf Proof.} We proceed similarly as in the proof of Lemma \ref{3momentep}. We need to establish that $$\int_{\Td} \rho(t,x)(\nabla\psi(t,x) -x)dx \equiv 0.$$
For this we use Definition \ref{3MaTdrho}. For $f(x)= \nabla\psi(x)-x$, we have
\be
\int_{\Td} \rho f& =& \int_{\Td} f(\nabla\hat\phi)\\
&=& \int_{\Td}\nabla\psi(\nabla\phi )-\nabla\phi \\
&=& \int_{\Td} x-\nabla\phi\\ 
&=&0,
\en
where we have used at the third line that, for $\psi, \phi$ Legendre transform of each other, $\nabla\psi(\nabla\phi) = id$, and at the last line that  $\phi -|x|^2/2$ is periodic. 

$\hfill \Box$

Hence we have shown that one can retrieve $v$ from the initial value of $\int \rho v$, and $\nabla\cdot v, \text{curl} v$.

\paragraph*{General framework}
We perform the same div-curl decomposition as in the Euler-Poisson case.
We then express the difference between the solution of $(EMA\epu)$ and the limiting solution : either the solution of $(E)$ or the solution of $(EMA\epu)$. After having applied a proper scaling to this difference, 
our solution is now described by a vector ${\bf u}$ whose first component (that can be a vector if $d=3$) is the rescaled vorticity, and whose last two components are a rescaled divergence and rescaled density fluctuation.
For this perturbation we will obtain 
\be
&&\dt {\bf u}\ep + \sum_i v^i \partial_i{\bf u}\ep + R\ep {\bf u}\ep = Q\ep({\bf u}\ep),
\en
where we still use
\be
R\ep=\left(\begin{array}{ccc}0&0&0\\0&0&-\frac{1}{\epsilon}\\0&\frac{1}{\epsilon}&0\end{array}\right).
\en
For the source term $Q\ep$, we have $\|Q\ep\| \leq C(1+\|{\bf u}\ep\|_{H^s}+\delta_\epsilon\|{\bf u}\ep\|_{H^s})$ 
where $\delta_\epsilon$ goes to 0 when $\epsilon$ goes to 0, and the constant $C$
depends on  the regularity of the limiting field. 

{\it Regularity of the limiting field} The form of the source term will vary under the circumstances, but  the general idea is that in order to bound $Q\ep$ in $H^s$, we will need the limiting velocity to be bounded in $H^{s+2}$ and the limiting density to be bounded in $H^{s+1}$. Remember that the $H^s$ norm of ${\bf u}$ controls 
the norm  of $(\rho, v)$ in $H^{s}\times H^{s+1}$, 
thus the limiting field must have one more derivative bounded than the order of the energy estimate. 

A Gronwall's lemma then yields a control on the perturbation that holds on a range of time $[0,T\epu]$, where $T\epu \to T$, $T$ being the time of existence of a smooth solution for  the limiting equation.

\paragraph*{Convergence to Euler, two dimensional case}
Doing the same change of variables as in the proof of Theorem \ref{3eulerpoisson} 
\be
&& \nabla \cdot v = \epsilon \beta_1,\\
&& \rho=1+ \epsilon^2 \rho_1, \\
&& \textrm{curl} v =\omega=\bar \omega+\epsilon \omega_1,
\en
we obtain:
\beq
&&\dt (\bar \omega+\epsilon \omega_1)+ v\cdot \nabla (\bar \omega+\epsilon \omega_1)
=-(\bar \omega+\epsilon \omega_1)\epsilon \beta,\label{3vortimaeu}\\
&& \dt \epsilon\beta_1 + v\cdot\nabla \epsilon\beta_1 + 2\epsilon\partial_i \bar v^j \partial_j v_1^i +\epsilon\partial_i v_1^j \partial_j v_1^i=\frac{\Delta \psi -d}{\epsilon^2 }-\partial_i \bar v^j \partial_j \bar v^i,\\
&& \dt \epsilon^2\rho_1 + v\cdot \nabla\epsilon^2\rho_1 = -(1+\epsilon^2\rho_1)\epsilon \beta_1.
\enq
Now   we define $\Xi$ by 
\be
\Delta\psi-d =\epsilon^2\rho_1+ \epsilon^4\Xi,
\en 
and from Theorem \ref{3maHs} inequality (\ref{3detdelta2}), 
we have, if $s\geq 2$, 
$\|\Xi\|_{H^s}\leq C \|\rho_1\|^2_{H^s}$.
The system can here be written in the following way:
\be
&&\dt {\bf u}\ep + \sum_i v^i \partial_i{\bf u}\ep + R\ep {\bf u}\ep = S\ep({\bf u}\ep)
+V\ep,\\
&& {\bf u}\ep(0)={\bf u}\ep_0,
\en
still with 
\be
{\bf u}\ep=\left(\begin{array}{c}\omega_1\\\beta_1\\ \tilde \rho_1 \end{array}\right)
, R\ep=\left(\begin{array}{ccc}0&0&0\\0&0&-\frac{1}{\epsilon}\\0&\frac{1}{\epsilon}&0\end{array}\right),
\en
with the same $S\ep$ as in the Euler-Poisson case, $\tilde \rho_1=\rho_1-\Delta p$, and with 
\be
V\ep=\left(\begin{array}{ccc}0\\\epsilon \Xi \\0\end{array}\right).
\en
We have $\|V\ep\|_{H^s}\leq C\epsilon (1+ \|{\bf u}\ep\|^2_{H^s})$, for $s$ large enough.
Then the energy estimates are the same as in the first proof, the solution 
${\bf u}\ep$
satisfying a control of the form:
\be
\Dt \|{\bf u}\ep\|^2_{H^s}\leq C \left(1+ \|{\bf u}\ep\|^2_{H^s} + \epsilon \|{\bf u}\ep\|^3_{H^s}\right)
\en
and the same conclusion holds true.
Then from Corollary \ref{3estime}, $v-\bar v, \rho-1$ can be retrieved from ${\bf u}\ep$, and we obtain the expected conclusion.

$\hfill\Box$

\paragraph*{Convergence to Euler, three dimensional case}

In the 3-d case, equation (\ref{3vortimaeu}) should be replaced by 
\be
\dt (\bar \omega+\epsilon \omega_1)+ v\cdot \nabla (\bar \omega+\epsilon \omega_1)
=-(\bar \omega+\epsilon \omega_1)\epsilon \beta
-\omega\cdot \nabla v.
\en
Note that the vorticity equation is the same as in the Euler-Poisson case.
This change would not affect the energy estimates.

\paragraph*{Higher order approximation}

Here we prove that the the Euler-Poisson system and the Euler-Monge-Amp\`ere system are closer as $\epsilon$ goes to 0 than Euler-Poisson and Euler.
We fix $s\geq s_0=[d/2]+2$.
For $\bar v_0$ a $H^{s+2}$ smooth divergence-free vector field, we consider $(\bar v, p)$  a solution of the Euler incompressible system (\ref{3euler}) such that $\bar v \in \Linf([0,T], H^{s+2}(\Td))$ for some $T>0$.
We consider also a sequence $(v\ep_{ep},\rho\ep_{ep})$ of solutions of the $(EP\epu)$ system with initial data $(v\ep_{ep,0},\rho\ep_{ep,0})$ (with $\int \rho_0\ep=1$) such that  $\epsilon^{-1}(v\ep_{ep}-\bar v),\epsilon^{-2}(\rho\ep_{ep}-1)$
is bounded in $\Linf([0,T'], H^{s+1}\times H^{s}(\Td))$, for any  $0<T'<T$,  if $\epsilon$ is small enough.
Thanks to Theorem \ref{3eulerpoisson}, and from the regularity assumption made on $\bar v$,  such a sequence exists for any sequence of well prepared initial data.

\begin{theo}\label{3poisson-ma}
Let $s\in\N$ with $s\geq [d/2]+2$. Let $\bar v,v\ep_{ep},\rho\ep_{ep}$ be as above.
Let $ (v_0\ep, \rho_0\ep)$ be a sequence of initial data such that
$\ds\left(\epsilon^{-2}(v_0\ep-v\ep_{ep,0}),\epsilon^{-3}(\rho_0\ep-\rho\ep_{ep,0})\right)$ 
is bounded in $H^s\times H^{s-1}(\Td)$. Then there exists a sequence $ (v\ep, \rho\ep)$ of 
solutions to $(EMA\epu)$  with initial data $(v_0\ep, \rho_0\ep)$, belonging to  $\Linf([0,T_{\epsilon}],H^s\times H^{s-1}(\Td))$, 
with $\liminf_{\epsilon\rightarrow 0}T_{\epsilon}\geq T$. Moreover for any $T'<T$ and $\epsilon$  small enough, the sequence 
$\ds\left(\epsilon^{-2}(v\ep-v\ep_{ep}), \epsilon^{-3}(\rho\ep-\rho\ep_{ep})\right)$
is bounded in $\Linf([0,T'], H^{s}\times H^{s-1}(\Td))$.
\end{theo}

{\it Remark.} We see here that the difference between solutions of $(EP\epu)$ and $(EMA\epu)$ is of order
$\epsilon^3$ for the density and of order $\epsilon^2$ for the velocity whereas the difference
between solutions of $(EP\epu)$ (or $(EMA\epu)$) and Euler was of order $\epsilon^2$ for the density and of order $\epsilon$ for the velocity.

\bigskip

{\bf Proof.}
We give the proof when $d=2$, the proof would be the same when $d=3$, just with more terms.
We introduce $(v_{ep},\rho_{ep}=1+\epsilon^2\rho_1)$ solution to $(EP\epu)$, and $(\beta_{ep}, \omega_{ep})=(\nabla\cdot v_{ep}, \nabla\times v_{ep})$.
Then we set
\be
&& v=v_{ep} + \epsilon^2 v_2,\\
&& \nabla \cdot v = \beta_{ep}+ \epsilon^{2}\beta_2,\\
&& \rho=\rho_{ep} + \epsilon^{3}\rho_2,\\
&& \textrm{curl} v =\omega=\omega_{ep} + \epsilon^2 \omega_2. 
\en
The system \( (EMA\epu) \) now reads:
\beq
&&\dt (\omega_{ep}+\epsilon^2 \omega_2 )+ v\cdot \nabla (\omega_{ep}+\epsilon^2\omega_2 )
=-(\omega_{ep}+\epsilon^2 \omega_2 )(\beta_{ep} + \epsilon^2\beta_2),\\
&& \dt (\beta_{ep} + \epsilon^{2}\beta_2)+ v\cdot\nabla (\beta_{ep}+\epsilon^{2}\beta_2)
+ \nabla(v_{ep}+ \epsilon^2 v_2):\nabla(v_{ep}+ \epsilon^2 v_2)
=\frac{\Delta \psi -d}{\epsilon^2 },\\
&& \dt (\rho_{ep} + \epsilon^{3}\rho_2)+ v\cdot \nabla(\rho_{ep}+ \epsilon^{3}\rho_2) 
= -(\rho_{ep}+ \epsilon^{3}\rho_2)(\beta_{ep}+ \epsilon^{2}\beta_2).
\enq
We still define 
$\Xi$ by 
\be
\Delta\psi-d =\rho-1 + \epsilon^4\Xi,
\en 
and from Theorem \ref{3maHs} we will have
\be
\|\Xi\|_{H^s(\Td)}\leq C \epsilon^{-4}\|\rho-1\|^2_{H^s(\Td)} \leq C(\|\rho_1\|^2_{H^s(\Td)} + \epsilon^2\|\rho_2\|^2_{H^s(\Td)}),
\en
(we use the notation $\rho_{ep}=1+ \epsilon^2 \rho_1$).

Setting
\be
{\bf u}\ep=\left(\begin{array}{c}\omega_2\\\beta_2\\ \rho_2 \end{array}\right),
\en
 we obtain that 
\be
\dt {\bf u}\ep + v\cdot \nabla {\bf u}\ep + R\ep {\bf u}\ep = T\epu,
\en
with $R\ep$ as before and $T\epu$ defined by
\be
T\epu=\left(\begin{array}{c}  -v_2\cdot \nabla w_{ep} -\beta_{ep} \omega_2-\beta_2 \omega_{ep} -\epsilon^2\omega_2 \beta_2\\
-v_2\cdot \nabla \beta_{ep} -2 \nabla v_{ep} : \nabla v_2 -\epsilon^2 \nabla v_2 : \nabla v_2 + \Xi\\
-\epsilon v_2\cdot \nabla \rho_1 -\beta_{ep}\rho_2 -\epsilon\beta_2\rho_1 -\epsilon^2\beta_2\rho_2
\end{array}\right). 
\en

Using again Proposition \ref{3AG} as in Lemma \ref{3bound} we obtain that, for $s>d/2$,
\be
\|T\epu\|_{H^s(\Td)}\leq C_s(1+\|{\bf u}\ep\|_{H^s(\Td)}+ \epsilon\|{\bf u}\ep\|_{H^s(\Td)}^2),
\en 
where the constant $C_s$ is controlled by $\|v_{ep}\|_{H^{s+2}}, \|\rho_1\|_{H^{s+1}}$ (still with $\rho_1=\epsilon^{-2}(\rho_{ep}-1)$).  From Theorem \ref{3eulerpoisson} these quantities are controlled 
for $0\leq t\leq T'<T$,  $T$ being  the time on which the solution of $(E)$ is smooth.  Hence we have by Gronwall's lemma a bound on $\|{\bf u}\ep\|_{\Linf([0,T'], H^s)}$.

Arguing as in the previous proofs, and using Corollary \ref{3estime}, we obtain that $( v_2, \rho_2)$ remains bounded in $\Linf([0,T''], H^{s+1}\times H^{s}(\Td))$ for any $T''<T'$ and for  $\epsilon<\epsilon_0$ small enough.
It follows that $(\epsilon^{-2}(v_{ep}-v_{ema}), \epsilon^{-3}(\rho_{ep}-\rho_{ema}))$
remains bounded in $\Linf([0,T''], H^s\times H^{s-1}(\Td))$. This achieves the proof of Theorem \ref{3poisson-ma}.

$\hfill \Box$

\subsection{Non-prepared initial data}
In this case, we obtain exactly the same result as for Euler-Poisson, using the same techniques. We follow closely section 2.3, and 
we only have to estimate the additional source term that will appear in the equation followed by $\beta_1$, due to the Monge-Amp\`ere coupling.
We recall that $\det D^2\psi=\rho$, and
we will have to estimate the difference
\be
\frac{1}{\epsilon^2}\left[( \Delta \psi - d) - (\rho-1)\right]
\en
in $H^s$ when we know that $\epsilon^{-1} (\rho-1)$ is bounded in $H^s$.
Thanks to Theorem \ref{3maHs}, we conclude that for $s\geq 2$ this term is controlled by $\|\epsilon^{-1}(\rho-1)\|^2_{H^s} = \|\rho_1\|_{H^s}$. Hence the energy estimate can be handled similarly just with an additional term, and the conclusion remains true.

\begin{theo}\label{3mongeamperenp}
The Theorem \ref{3eulerpoissonnp} holds also when replacing the $(EP\epu)$ system by the $(EMA\epu)$ system.
\end{theo}

\( \hfill \Box \)

\bibliography{ema-biblio}
\vspace{1cm}
\begin{flushright} 
Gr\'egoire Loeper
\\
EPFL, SB-IMA
\\ 
1015 Lausanne, Switzerland
\\
gregoire.loeper@epfl.ch
\end{flushright}

\end{document}